\documentclass[12pt]{amsart}
\usepackage{amsmath,amscd,amssymb,amsfonts}
\newcommand{\h}{\hbox}
\newcommand{\q}{\quad}

\newcommand{\bs}{\par\bigskip}
\newcommand{\ms}{\par\medskip}
\newcommand{\sk}{\par\smallskip}
\newcommand{\bsn}{\par\bigskip\noindent}
\newcommand{\msn}{\par\medskip\noindent}
\newcommand{\skn}{\par\smallskip\noindent}
\newcommand{\msum}{\h{$\sum$}}
\newcommand{\mopl}{\h{$\bigoplus$}}
\newcommand{\mcup}{\h{$\bigcup$}}
\newcommand{\mprod}{\h{$\prod$}}
\newcommand{\ssb}{\raise.15ex\h{${\scriptscriptstyle\bullet}$}}
\newcommand{\ssc}{\,\raise.15ex\h{${\scriptstyle\circ}$}\,}
\newcommand{\ds}{\rlap{\raise.4ex\h{$\downarrow$}}{\downarrow}}
\newcommand{\C}{{\mathbf C}}
\newcommand{\N}{{\mathbf N}}
\newcommand{\PP}{{\mathbf P}}
\newcommand{\Q}{{\mathbf Q}}

\newcommand{\Z}{{\mathbf Z}}
\newcommand{\D}{{\mathcal D}}
\newcommand{\F}{{\mathcal F}}
\newcommand{\Hc}{{\mathcal H}}
\newcommand{\OO}{{\mathcal O}}
\newcommand{\mm}{{\mathfrak m}}
\newcommand{\Ht}{\widetilde{H}}
\newcommand{\Yt}{\widetilde{Y}}
\newcommand{\Zt}{\widetilde{Z}}
\newcommand{\dd}{{\partial}}
\newcommand{\ddd}{{\mathrm d}}
\newcommand{\dddd}{\overline{\rm d}}
\newcommand{\df}{{\mathrm d}f}
\newcommand{\al}{\alpha}
\newcommand{\la}{\lambda}
\newcommand{\om}{\omega}
\newcommand{\Om}{\Omega}
\newcommand{\sw}{\,\h{$\wedge$}\,}
\newcommand{\EV}{{\rm EV}}
\newcommand{\Gr}{{\rm Gr}}
\newcommand{\Sp}{{\rm Sp}}
\newcommand{\bl}{\bigl}
\newcommand{\br}{\bigr}
\newcommand{\into}{\hookrightarrow}
\newcommand{\onto}{\mathop{\rlap{$\to$}\hskip2pt\hbox{$\to$}}}
\newcommand{\simto}{\buildrel\sim\over\too}
\newcommand{\too}{\longrightarrow}
\newcommand{\ges}{\geqslant}
\newcommand{\les}{\leqslant}
\begin{document}
\title[Hilbert series of Milnor algebras]
{Hilbert series of graded Milnor algebras and roots of Bernstein-Sato polynomials}
\author[M. Saito]{Morihiko Saito}
\begin{abstract}
We show that there is a pair of homogeneous polynomials such that the sets of roots of their Bernstein-Sato polynomials which are strictly supported at the origin are different although the sets of roots which are not strictly supported at the origin are the same and moreover their graded Milnor algebras have the same Hilbert series. This shows that the roots of the Bernstein-Sato polynomials strictly supported at the origin cannot be determined uniquely by the Hilbert series of the Milnor algebras. This is contrary to certain hyperplane arrangement cases. It also implies that a nonzero torsion element with pure degree in the Milnor algebra does not necessarily contribute to a root of the Bernstein-Sato polynomial in an expected way. This example is found by using Macaulay2 and RISA/ASIR.
\end{abstract}
\maketitle
\centerline{\bf Introduction}
\bsn
Let $f$ be a homogeneous polynomial of $n$ variables with degree $d$. In this paper we assume
$$\h{$Z:=\{f=0\}\subset\PP^{n-1}$ has only isolated singularities.}
\leqno(1)$$
Let $\Om^{\ssb}$ be the graded complex of algebraic differential forms on $X:=\C^n$.
Set
$$M:=H^n(\Om^{\ssb},\df\wedge),\q N:=H^{n-1}(\Om^{\ssb},\df\wedge),\q M':=H^0_{\mm}M,\q M''=M/M',$$
where ${\mm}:=(x_1,\dots,x_n)\subset\C[x]:=\C[x_1,\dots,x_n]$ with $x_1,\dots,x_n$ the coordinates of $\C^n$, and the grading of $\Om^j$ is shifted by $d(n-j)$ so that $\df\wedge$ preserves the grading without a shift, see \cite{DS2}. (In fact, $(\Om^{\ssb},\df\wedge)$ is identified with a graded piece of a filtered Gauss-Manin complex.)
Note that $M$ is identified with the {\it Milnor algebra} $\,\C[x]/(\dd f)$ up to the shift of degree by $n$, where $(\dd f)$ is the Jacobian ideal generated by the partial derivatives of $f$.
In this paper, $M'$ is called the {\it torsion part} of $M$. This can be justified by considering $M$ over $\C[y]$ with $y$ a sufficiently general linear combination of the coordinates $x_i$ of $\C^n$.
\sk
Let $\EV(Z):=\mcup_{z\in{\rm Sing}\,Z}\,\EV(Z,z)$ with $\EV(Z,z)$ the set of the eigenvalues of the Milnor monodromy of a local equation $h_z$ of $Z$ at $z\in{\rm Sing}\,Z$. Set
$$R^{\,0}_f:=\{\al\in\Q\mid b_f(-\al)=0,\,\,\exp(-2\pi i\al)\notin\EV(Z)\}\subset\Q_{>0},$$
(see \cite{Ka1} for the last inclusion). This is called the set of roots of $b_f(-s)$ which are {\it strictly supported at the origin}.
The following may be viewed as a partial generalization of a result of Malgrange \cite{Ma1} in the isolated singularity case.
\msn
{\bf Theorem~1} (\cite[Thm.~2]{bCM}). {\it Let $k$ be a positive integer. Assume
$$\exp(-2\pi ik/d)\notin\EV(Z).
\leqno(2)$$
Then
$$k/d\in R^{\,0}_f\iff\Gr^p_PH^{n-1}(f^{-1}(1),\C)_{\la}\ne 0,
\leqno(3)$$
where $\la=\exp(-2\pi ik/d)$, $p=[n-k/d]$, and $P$ is the pole order filtration on the $\la$-eigenspace of the monodromy on the Milnor cohomology.
We have moreover a canonical surjection
$$M_k\onto M_k^{(\infty)}=\Gr^p_PH^{n-1}(f^{-1}(1),\C)_{\la},
\leqno(4)$$
with $\la$, $p$ as above by using the pole order spectral sequence, see $(2.1)$ below.}
\ms
Note, however, that it is {\it quite difficult} to determine the pole order filtration $P$ on the Milnor cohomology by using only the pole order spectral sequence without using a computer program except for certain simple cases as in \cite{bha}, see (2.4) below.
\sk
By (3), (4) we say that an element of $M_k$ {\it contributes} to a root $-k/d$ of the Bernstein-Sato polynomial $b_f(s)$ if its image in $M_k^{(\infty)}$ does not vanish. Here we also assume condition~(2).
\sk
In this paper, we show the following
\msn
{\bf Theorem~2.} {\it The set $R^{\,0}_f$ is not uniquely determined by the Hilbert series $\{\mu_{k+n}\}_{k\in\N}$ of the Milnor algebra $\C[x]/(\dd f)$.}
\ms
This implies that the differentials $\ddd_r$ of the pole order spectral sequence are not uniquely determined by the $E_1$-term of the spectral sequence.
For the proof of Theorem~2, we study a pair of homogeneous polynomials
$$f_1=x^5+y^4z+x^4y,\q\q f_2=x^5+y^4z+x^3y^2,$$
where $n=3$, $d=5$. (We have been informed that $f_1$ already appeared in \cite[Remark 2.8(5)]{Wa}.)
These have the same Hilbert series of their Milnor algebras according to computer calculations using Macaulay2. On the other hand, calculations using RISA/ASIR imply
$$8/5\in R^{\,0}_{f_2}\setminus R^{\,0}_{f_1},
\leqno(5)$$
see (1.4) below. So Theorem~2 follows.
As another consequence of (5), we get the following.
\msn
{\bf Theorem~3.} {\it For certain homogeneous polynomials $f$ and integers $k$ in $[d+1,2d-2]$ satisfying conditions~$(1)$ and $(2)$, it is possible that we have $k/d\notin R^{\,0}_f$ although $M'_k\ne 0$ and $M'_{k+1}=0$.}
\ms
Note that the $\dim M'_k$ ($k\in\N$) are uniquely determined by the $\dim M_k$ ($k\in\N$), see (2.4.1) below. In the cases of $f_1$ and $f_2$, we have
$$\dim M'_8=1,\q\dim M'_9=0.
\leqno(6)$$
This is closely related to an assertion in an old version of \cite{Wa}, since $M'_8$ is annihilated by the maximal ideal $(x,y,z)$ by the vanishing of $M'_9$. Here it should be noted that there is no canonical morphism
$$M'\to H^{n-1}(f^{-1}(1),\C),
\leqno(7)$$
see also (10) below. In fact, we have no canonical isomorphism even in the isolated singularity case where $M=M'$ and a canonical isomorphism (7) would imply a canonical opposite filtration to the Hodge filtration, see for instance \cite{bl}.
\sk
Related to Theorem~3, the following is proved in this paper without using a computer.
\msn
{\bf Theorem~4.} {\it For a certain homogeneous polynomial $f$ satisfying condition~$(1)$, there is a homogeneous polynomial $g$ of degree $k-n$ such that $k$ is contained in $[d+1,2d-2]$, and satisfies condition~$(2)$, and we have
$$M'_k\ni[g\ddd x]\ne 0,\q\h{but}\q[g\,\widehat{\ddd x}|_{f=1}]=0\q\h{in}\,\,\,H^{n-1}(f^{-1}(1),\C),
\leqno(8)$$
where $\ddd x:=\ddd x_1\sw\cdots\sw\ddd x_n$, and $\widehat{\ddd x}$ is the contraction of $\ddd x$ with the Euler field $\msum_{i=1}^n\,x_i\dd_{x_i}$. Moreover $[g\ddd x]\in M'_k$ belongs to the kernel of the canonical morphism
$$M'_k\into M_k\onto M_k^{(\infty)}=\Gr^p_PH^{n-1}(f^{-1}(1),\C)_{\la},
\leqno(9)$$
where $\la$, $p$ are as in $(3)$.}
\ms
Note that a similar assertion holds with condition $(8)$ replaced by
$$[g\ddd x]=0\q\h{in}\,\,\,M'_k,\q\h{but}\q[g\,\widehat{\ddd x}|_{f=1}]\ne 0\q\h{in}\,\,\,H^{n-1}(f^{-1}(1),\C),
\leqno(10)$$
see the last remark in (2.3) below.
In fact, this may happen even in the isolated singularity case, and is closely related with the ambiguity of the morphism (7).
\ms
Theorem~4 means that each nonzero element of $M'_k$ does not necessarily {\it contribute} to a root $-k/d$ of $b_f(s)$.
(Note, however, that ``does not contribute to a root $-k/d\,$" does not necessarily mean that $-k/d$ is not a root of $b_f(s)$.)
\sk
For the proof of Theorem~4, we study the case of the the above polynomial $f_1$ with $k=8$. We calculate part of the pole order spectral sequence (see (2.1) below) for $f_1$, and show the inclusion
$$M'_8\subset{\rm Im}(\ddd:N_{13}\to M_8),
\leqno(11)$$
which is closely related to the last assertion of Theorem~4. More precisely, set
$$\om:=x^2y^3\,\om_0\q\h{with}\q\om_0:=\ddd x\sw\ddd y\sw\ddd z.$$
By (6), the inclusion (11) is reduced to
$$M'_8\cap{\rm Im}(\ddd:N_{13}\to M_8)\ni[\om]_M\ne 0,
\leqno(12)$$
see (2.3) below. Here $[\om]_M$ denotes the class of $\om$ in $M$.
Note that (12) implies the vanishing of the image of $[\om]_M$ in $M_8^{(\infty)}$, and the last assertion of Theorem~4 follows, where $g$ is defined as in (2.3.6) below for the first assertion of Theorem~4 (by using that $x^3y^2\in(\dd f)$).
\sk
Note that $[\om]_M$ contributes to the root $-3/5$ of $b_f(s)$ in some sense, since
$$[\om]_{H''_f}=c\,\dd_t^{-1}[\om_0]_{H''_f}\q(c\in\C^*).$$
Here $[\om]_{H''_f}$ denotes the class of $\om$ in the {\it algebraic Brieskorn module} $H''_f:=H^3(A_f^{\ssb},\ddd)$ in the notation of (2.2) below (see \cite{Br}, \cite{BaS}, etc.\ for the analytic case). We have
$$\h{$\frac{3}{5}$}\,\dd_t^{-1}[\om_0]_{H''_f}=t\,[\om_0]_{H''_f},$$
and $[\om_0]_M$ contributes to the root $-3/5$ of $b_f(s)$. In fact, we know that $-3/5$ is a root by a computer calculation although this does not hold for $f_2$, see (1.4) below.
\sk
We thank A.~Dimca and U.~Walther for their useful comments and for calculating the Hilbert series of the Milnor algebras and the Bernstein-Sato polynomials for certain examples treated in an earlier version of this paper, see Remark~(2.5) below.
This work is partially supported by Kakenhi 24540039.
\sk
In Section~1 we calculate some numerical invariants of the singularities of $Z$ for $f_1$.
In Section~2 we prove Theorem~4, and then calculate the pole order spectral sequence for $f_1$ and $f_2$.
In Appendix we give a proof of a variant of an assertion in an old version of \cite{Wa} without assuming the knowledge of logarithmic differential forms very much for nonspecialists.
\bs\bs
\vbox{\centerline{\bf 1. Preliminaries}
\bsn
In this section we calculate some numerical invariants of the singularities of $Z$ for $f_1$. To simplify the notation, we will denote $f_1$ by $f$ in this section except in (1.4).}
\msn
{\bf 1.1.~Singularities of $Z$.}
We first see that $Z$ is a rational curve having only one singular point $z_0$ at $[0:0:1]$. Indeed, substituting respectively $x=1$ and $y=1$ to $f=0$, we get
$$(1+y)+y^4z=0\q\h{and}\q x^4(x+1)+z=0.$$
These define nonsingular rational curves in $\C^2$. Set
$$h:=x^5+x^4y+y^4.$$
At the singular point $z_0$ of $Z$, we have
$$\mu_h=12,\q\tau_h=11.
\leqno(1.1.1)$$
\msn
{\bf 1.2. Some numerical invariants.}
Since the singular point $z_0$ of $Z$ is analytic-locally irreducible, we have
$$H^j(Z,\C)=\begin{cases}\C&\h{if}\,\,\,\,j=0,2,\\ \,0&\h{if}\,\,\,\,j\ne 0,\,2.\end{cases}
\leqno(1.2.1)$$
Set
$$U:=\PP^2\setminus Z.$$
By using the long exact sequence associated with the cohomology with compact support together with duality, we get
$$\widetilde{H}^j(U,\C)=0,
\leqno(1.2.2)$$
\sk
Let $H^j(f^{-1}(1),\C)_{\la}$ be the $\la$-eigenspace of the monodromy on the Milnor cohomology. We have
$$H^1(f^{-1}(1),\C)_{\la}=0\q(\forall\,\la),
\leqno(1.2.3)$$
which is equivalent to
$$\nu_k^{(\infty)}=0\q(\forall\,k).
\leqno(1.2.4)$$
In fact, (1.2.3) for $\la=1$ follows from (1.2.2). For $\la\ne 1$, it is enough to show that, for any eigenvalue $\la'\ne 1$ of the Milnor monodromy of the above $h$, we have $\la'{}^{\,5}\ne 1$. Here we can replace $h$ with a Brieskorn-Pham polynomial $h':=x^5+y^4$ by using a $\mu$-constant deformation. Then the assertion is well-known (see also the calculation of the Steenbrink spectrum of $h'$ in (1.3) below).
\sk
In the case of a homogeneous polynomial of degree $d$, it is well-known that
$$H^j(f^{-1}(1),\C)_{\la}=0\q\h{unless}\q\la^d=1.$$
More precisely there are local systems $L_{\la}$ of rank 1 on $U$ for $\la^d=1$ such that $L_1=\C_U$ and
$$H^j(f^{-1}(1),\C)_{\la}=H^j(U,L_{\la}),$$
see for instance \cite{Di}. The above calculations then imply
$$\dim H^2(f^{-1}(1),\C)_{\la}=\begin{cases}1&\h{if}\,\,\,\,\la^5=1,\,\la\ne 1,\\ 0&\h{otherwise.}\end{cases}
\leqno(1.2.5)$$
\msn
{\bf 1.3.~Spectrum.} We first recall the definition of spectrum in the general case. For a holomorphic function $f$ on a germ of a complex manifold $(X,0)$ of dimension $n$, we have the $j$\,th Steenbrink spectrum $\Sp^j(f)=\msum_{\al\in\Q}\,n^j_{f,\al}\,t^{\al}$ defined by
$$n^j_{f,\al}:=\dim\Gr_F^p\Ht^{n-1-j}(F_{\!f,0},\C)_{\la}\q\bl(p:=[n-\al],\,\la=\exp(-2\pi i\al)\br),$$
where $F$ is the Hodge filtration on the $\la$-eigenspace of the monodromy of the reduced Milnor cohomology $\Ht^j(F_{\!f,0},\C)_{\la}$ with $F_{\!f,0}$ the Milnor fiber of $f$ around 0, see \cite{St2}, \cite{St3}.
The Steenbrink spectrum $\Sp^j(f)$ is given by
$$\Sp(f)=\msum_{j=0}^{n-1}\,(-1)^j\,\Sp^j(f).$$
Similarly we have the $j$\,th pole order spectrum
$$\Sp^j_P(f)=\msum_{\al\in\Q}\,n^j_{P,f,\al}\,t^{\al},$$
defined by replacing the Hodge filtration $F$ with the pole order filtration $P$, which satisfies
$$F^p\subset P^p.
\leqno(1.3.1)$$
This implies for instance
$$n_{P,f,\al}=0\q\h{if}\,\,\,\,n_{f,\al+i}=0\,\,\,\,\h{for any}\,\,\,\,i\in\N.
\leqno(1.3.2)$$
\sk
In the case of isolated singularities (or more generally, if $\Ht^j(F_{\!f,0},\C)=0$ for $j<n-1$), we have $\Sp(f)=\Sp^0(f)$, and similarly for $\Sp_P(f)$.
The rational numbers $\al$ with $n_{f,\al}\ne 0$ are called the Steenbrink spectral numbers. These are counted with multiplicities given by $n_{f,\al}$.
\sk
In our case, we can replace $h=x^5+x^4y+y^4$ with $h'=x^5+y^4$ for the calculation of the Steenbrink spectrum, since $h$ is a $\mu$-constant deformation of $h'$, see for instance \cite{Va1}. By \cite{St1}, \cite{St2}, the Steenbrink spectral numbers of $h'$ (and $h$) are given by
$$\aligned&\q\q\bl\{(4i + 5j)/20\mid (i,j)\in[1,4]\times[1,3]\}\br\}\\
&=\bl\{9,13,14,17,18,19,21,22,23,26,27,31\br\}/20.\endaligned
\leqno(1.3.3)$$
These do not intersect ${\bf Z}\frac{1}{5}\setminus{\bf Z}$.
This implies (1.2.3) since $H^1(U,\C)=0$ by (1.2.2).
\sk
Let $V$ be the $V$-filtration of Kashiwara \cite{Ka2} and Malgrange \cite{Ma2} on
$$(i_f)^{\!\D}_*\OO_X=\OO_X[\dd_t],$$
where the latter is the direct image of the structure sheaf $\OO_X$ as a left $\D$-module by the graph embedding $i_f$ of $f:X\to\C$ (and the sheaf-theoretic direct image on the right-hand side is omitted to simplify the notation).
\sk
Let $f=x^5+x^4y+y^4z$ as in the introduction. We have the induced filtration $V$ on
$$\C\{x,y,z\}=\OO_{X,0}\subset(i_f)^{\!\D}_*\OO_{X,0}=\OO_{X,0}[\dd_t].$$
By \cite[Thm.~2.2]{bCM} together with the relation with the multiplicative ideals (see \cite{Bu}, \cite{BS}), we see that $V^{\al}\C\{x,y,z\}$ is generated by $x^iy^jz^k$ satisfying the two conditions:
$$(4i+5j+9)/20\ges\al,\q(i+j+k+3)/5\ges\al.$$
For $\al=3/5$ and $4/5$, we have respectively
$$\aligned(4i+5j+9)/20\ges 12/20\,\,&\Longrightarrow\,\,(i+j+k+3)/5>3/5,\\
(4i+5j+9)/20\ges 16/20\,\,&\Longrightarrow\,\,(i+j+k+3)/5>4/5.\endaligned$$
These imply that
$$\Gr_V^{\al}\OO_X=0\q\h{for}\,\,\,\al=3/5,\,\,4/5.$$
Hence $3/5$ and $4/5$ are not Steenbrink spectral numbers of $f$, see \cite{Bu}, \cite{BS}.
\sk
The multiplicity of $k/5$ as a Steenbrink spectral number of $f$ is then $1$ for $k=6,7,8,9$ (by using the symmetry of the direct factors of $\Gr^W_{\ssb}\psi_f\Q[n-1]$ supported on the origin as mixed Hodge modules).
\sk
Combining these with (1.2.2), we thus get
$$\Sp(f)=t^{6/5}+t^{7/5}+t^{8/5}+t^{9/5}.
\leqno(1.3.4)$$
\msn
{\bf 1.4.~Bernstein-Sato polynomials.} By using three different programs: bfct, bfunction, and ndbf.bf\underline{\h{ }}local in RISA/ASIR, we get the same answers as below:
$$\aligned b_{f_1}(s)&=\h{$\bl(s+\frac{9}{20}\br)\bl(s+\frac{11}{20}\br)\bl(s+\frac{13}{20}\br)\bl(s+\frac{14}{20}\br)\bl(s+\frac{17}{20}\br)\bl(s+\frac{18}{20}\br)\bl(s+\frac{19}{20}\br)\bl(s+\frac{21}{20}\br)$}\\ &\h{$\q\,\,\bl(s+\frac{22}{20}\br)\bl(s+\frac{23}{20}\br)\bl(s+\frac{26}{20}\br)\bl(s+\frac{27}{20}\br)\bl(s+\frac{3}{5}\br)\bl(s+\frac{4}{5}\br)\bl(s+\frac{6}{5}\br)\bl(s+\frac{7}{5}\br)(s+1)$,}\endaligned
\leqno(1.4.1)$$
$$\aligned b_{f_2}(s)&=\h{$\bl(s+\frac{9}{20}\br)\bl(s+\frac{11}{20}\br)\bl(s+\frac{13}{20}\br)\bl(s+\frac{14}{20}\br)\bl(s+\frac{17}{20}\br)\bl(s+\frac{18}{20}\br)\bl(s+\frac{19}{20}\br)\bl(s+\frac{21}{20}\br)$}\\ &\h{$\q\,\,\bl(s+\frac{22}{20}\br)\bl(s+\frac{23}{20}\br)\bl(s+\frac{26}{20}\br)\bl(s+\frac{27}{20}\br)\bl(s+\frac{4}{5}\br)\bl(s+\frac{6}{5}\br)\bl(s+\frac{7}{5}\br)\bl(s+\frac{8}{5}\br)(s+1)$,}\endaligned
\leqno(1.4.2)$$
$$\aligned b_h(s)&=\h{$\bl(s+\frac{9}{20}\br)\bl(s+\frac{11}{20}\br)\bl(s+\frac{13}{20}\br)\bl(s+\frac{14}{20}\br)\bl(s+\frac{17}{20}\br)\bl(s+\frac{18}{20}\br)\bl(s+\frac{19}{20}\br)\bl(s+\frac{21}{20}\br)$}\\ &\h{$\q\,\,\bl(s+\frac{22}{20}\br)\bl(s+\frac{23}{20}\br)\bl(s+\frac{26}{20}\br)\bl(s+\frac{27}{20}\br)(s+1)$,}\endaligned
\leqno(1.4.3)$$
where $f_1$, $f_2$ are as in the introduction, and $h$ is as in (1.1). Note that we have the equality
$$b_{h_1}(s)=b_{h_2}(s)\,(=b_h(s)),$$
by the theory of Brieskorn lattices \cite{bl} using \cite{Ma1}, where $h_a$ is obtained by restricting $f_a$ to $z=1$ for $a=1,2$.
(In fact, the minimal Steenbrink spectral number of $h$ is $9/11$, and $31/20$ is the only Steenbrink spectral number of $h$ which is bigger than $1+9/20$ by (1.3.3).)
We then get
$$\aligned&b_{f_1}(s)=b_h(s)\,\mprod_{i\in I_1}\,(s+i/5)\q\h{with}\q I_1:=\{3,4,6,7\},\\ &b_{f_2}(s)=b_h(s)\,\mprod_{i\in I_2}\,(s+i/5)\q\h{with}\q I_2:=\{4,6,7,8\}.
\endaligned
\leqno(1.4.4)$$
This is closely related with (3).
\sk
Combined with (1.3.4), these calculations imply
$$\aligned&\Sp_P(f_1)=t^{3/5}+t^{4/5}+t^{6/5}+t^{7/5},\\&\Sp_P(f_2)=t^{4/5}+t^{6/5}+t^{7/5}+t^{8/5}.\endaligned
\leqno(1.4.5)$$
Note that these are compatible with (1.3.2).
\bs\bs
\vbox{\centerline{\bf 2. Proof of Theorem~4}
\bsn
In this section we prove Theorem~4, and then calculate the pole order spectral sequence for $f_1$ and $f_2$. We first recall the notion of a pole order spectral sequence.}
\msn
{\bf 2.1.~Pole order spectral sequences.}
With the notation and the assumption of the introduction, we have inductively the morphisms
$$\ddd^{(r)}:N^{(r)}\to M^{(r)}\q(r\ges 1),$$
defined by the differentials $\ddd_r$ of the {\it pole order spectral sequence} (see \cite{DS2}). Here
$$\ddd^{(1)}:N^{(1)}\to M^{(1)}$$
coincides with
$$\ddd:N\to M,$$
and we set inductively
$$N^{(r)}:={\rm Ker}\,\ddd^{(r-1)},\q M^{(r)}:={\rm Coker}\,\ddd^{(r-1)}\q(r\ges 2).$$
Setting $M_k^{(\infty)}:=M_k^{(r)}$ ($r\gg 0$), etc., we have moreover the conditions:
$$\h{$\ddd^{(r)}$ preserves the degree up to the shift by $-rd$,}
\leqno(2.1.1)$$
$$M_k^{(\infty)}=\Gr^p_PH^{n-1}(f^{-1}(1),\C)_{\la},\q N_k^{(\infty)}=\Gr^p_PH^{n-2}(f^{-1}(1),\C)_{\la},
\leqno(2.1.2)$$
with $\la=\exp(-2\pi ik/d)$, $p=[n-k/d]$.
Here $H^j(f^{-1}(1),\C)_{\la}$ denotes the $\la$-eigenspace of the monodromy on the Milnor cohomology, and $P$ is the pole order filtration (see also \cite{bCM}).
These are essentially equivalent to the {\it pole order spectral sequence} (see \cite{DS2}), and will be so called in this paper.
\msn
{\bf 2.2.~Brieskorn modules.} Following \cite{Br}, we can define $H_f$, $H'_f$, $H''_f$ by
$$H_f:=H^{n-1}\Om_{X/S}^{\ssb},\q H'_f:=H^{n-1}(\df\wedge\Om^{\ssb}),\q H''_f:=H^n(A_f^{\ssb}),$$
where $A_f^j:={\rm Ker}(\df\wedge:\Om^j\to\Om^{j+1})$, and $\Om_{X/S}^{\ssb}:=\Om^{\ssb}/\df\wedge\Om^{\ssb-1}$ in the notation of the introduction. There are short exact sequences of complexes
$$0\to\df\wedge\Om^{\ssb}[-1]\to\Om^{\ssb}\to\Om_{X/S}^{\ssb}\to 0,$$
$$0\to A_f^{\ssb}\to\Om^{\ssb}\to\df\wedge\Om^{\ssb}\to 0.$$
inducing the canonical isomorphisms (by assuming $n\ges 2$)
$$\dd:H_f\simto H'_f,\q\dd:H'_f\simto H''_f.$$
\sk
On the other hand there are canonical morphisms
$$\alpha:H_f\to H'_f,\q\beta:H'_f\to H''_f,$$
induced by
$$\Om_{X/S}^{\ssb}\buildrel{\df\wedge}\over\too\df\wedge\Om^{\ssb},\q\df\wedge\Om^{\ssb-1}\into A_f^{\ssb}.$$
It is known that the action of $\dd_t^{-1}$ on the Brieskorn modules coincides with $\beta\ssc\dd^{-1}$, etc. The action of $t$ is defined by the multiplication by $f$.
\sk
We then get
$$\beta\ssc\alpha\bl([g\widehat{\ddd x}]_{H_f}\br)=(\deg f)\,t\,[g\ddd x]_{H''_f}\q\h{for}\,\,\,\,g\in\C[x],
\leqno(2.2.1)$$
where $\widehat{\ddd x}$ is the contraction of the Euler field $\sum_ix_i\dd_{x_i}$ with $\ddd x:=\ddd x_1\sw\cdots\sw\ddd x_n$, see Theorem~4.
In fact, (2.2.1) easily follows from $\sum_ix_i\dd_{x_i}f=(\deg f)\,f$.
\msn
{\bf 2.3.~Proof of Theorem~4.} Since $f=f_1:=x^5+x^4y+y^4z$, we have
$$f_x=5x^4+4x^3y,\q f_y=x^4+4y^3z,\q f_z=y^4.$$
We then get inductively
$$y^4,\,\,x^4y,\,\,x^3y^2,\,\,x^5,\,\,xy^3z\in(\dd f),
\leqno(2.3.1)$$
by using respectively
$$f_z,\,\,\,f_y,\,\,\,f_x,\,\,\,f_x,\,\,\,f_y.$$
\sk
Set
$$\om:=x^2y^3\,\om_0\q\h{with}\q\om_0:=\ddd x\sw\ddd y\sw\ddd z.$$
Then (2.3.1) implies
$$[\om]_M=[x^2y^3\,\om_0]_M\in M'_8,
\leqno(2.3.2)$$
where $[\om]_M$ denotes the class of $\om$ in $M$.
\sk
We have
$$[\om]_M\ne 0,\q\h{that is,}\q x^2y^3\notin(\dd f),
\leqno(2.3.3)$$
since
$$x^2y^3\notin(x^3,y^4,y^3z)\supset(\dd f).$$
\sk
Set $\eta:=\eta_1+\eta_2$ with
$$\aligned\eta_1:=g_1\ddd x\q\h{for}\q &g_1:=x^3/y,\\
\eta_2:=g_2\ddd y\q\h{for}\q &g_2:=(1/5)\,x^3/y.\endaligned$$
In the notation of the introduction, we can verify
$$\xi:=\df\wedge\eta\in A_f^2\subset\Om^2,
\leqno(2.3.4)$$
and
$$\ddd\xi=-\df\wedge\ddd\eta=c_1\,x^3y^2\om_0+c_2\,x^2y^3\om_0\in\ddd A_f^2\q(c_1,c_2\in\C^*),
\leqno(2.3.5)$$
since
$$f_y(g_1)_z-f_z(g_1)_y=x^3y^2,\q f_x(g_2)_z-f_z(g_2)_x=-(3/5)x^2y^3.$$
Then (12) follows from (2.3.5) since $x^3y^2\in(\dd f)$ by (2.3.1). For the first assertion of Theorem~4, set
$$g:=c_1\,x^3y^2+c_2\,x^2y^3,
\leqno(2.3.6)$$
and use (2.2.1). (If we set $g=x^3y^2$, then (10) holds by the argument at the end of the introduction.) This finishes the proof of Theorem~4.
\msn
{\bf 2.4.~Calculation of the pole order spectral sequence.} Set
$$\mu_k=\dim M_k,\q \mu'_k=\dim M'_k,\q \mu''_k=\dim M''_k,\q \nu_k=\dim N_k,$$
and similarly for $\mu_k^{(r)}$, etc. Define $\gamma_k$ by
$$\msum_k\,\gamma_kt^k=\bl(t+\cdots+t^4\br)^3.$$
Using the calculation of the Bernstein-Sato polynomials in (1.4), we can show the following:
\skn
For $f=f_1$, we have
$$\begin{array}{rccccccccccccc}
k:&3 &4 &5 &6 &7 &8 &9 &10 &11 &12 &13 &14 &\cdots\\
\gamma_k:&1 &3 &6 &10 &12 &12 &10 &6 &3 &1 &\\
\mu'_k:& & & & &1 &1 & & & & & &\\
\mu''_k:&1 &3 &6 &10 &11 &11 &11 &11 &11 &11 &11 &11 &\cdots\\
\nu_k:& & & & & & &1 &5 &8 &10 &11 &11 &\cdots\\
\mu_k:&1 &3 &6 &10 &12 &12 &11 &11 &11 &11 &11 &11 &\cdots\\
\nu_{k+5}:& &1 &5 &8 &10 &11 &11 &11 &11 &11 &11 &11 &\cdots\\
\mu^{(2)}_k:&1 &2 &1 &2 &2 &1 &1 &1 &1 & 1 & 1 &1 &\cdots\\
\nu^{(2)}_{k+5}:& & & & & & & 1 & 1 & 1 & 1 & 1 &1 &\cdots\\
\mu^{(3)}_k:&1 &1 & &1 &1 & &\\
\end{array}$$
\skn
For $f=f_2$, there are two possibilities; either
$$\begin{array}{rccccccccccccc}
k:&3 &4 &5 &6 &7 &8 &9 &10 &11 &12 &13 &14 &\cdots\\
\mu^{(3)}_k:&1 &1 & &1 &1 &\underline{1} &\\
\nu^{(3)}_{k+5}:& & & & & & & & & & &\underline{1} &\\
\mu^{(4)}_k:&\underline{0} &1 & &1 &1 &\underline{1} &\\
\end{array}$$
or
$$\begin{array}{rccccccccccccc}
k:&3 &4 &5 &6 &7 &8 &9 &10 &11 &12 &13 &14 &\cdots\\
\mu^{(2)}_k:&1 &2 &1 &2 &2 &\underline{2} &1 &1 &1 &1 &1 &1 &\cdots\\
\nu^{(2)}_{k+5}:& & & & & &\underline{1} &1 &1 &1 &1 &1 &1 &\cdots\\
\mu^{(3)}_k:&\underline{0} &1 & &1 &1 &\underline{1} &\\
\end{array}$$
and we do not know yet which is true.
Here some part of the table is omitted if it is entirely same as in the case of $f_1$, and the part different from the case of $f_1$ is underlined.
\sk
Recall that
$$\mu_k-\nu_k=\gamma_k,\q\mu''_k+\nu_{nd-k}=11\,(=\tau)\q(\forall\,k),
\leqno(2.4.1)$$
(see for instance \cite{DS2}). So we can calculate $\nu_k$, $\mu''_k$, $\mu'_k$ for any $k$.
\sk
By the argument in \cite[Section 5.1]{DS2} together with the calculation of the Steenbrink spectrum of $h$ in (1.3), we see that
$${\rm rank}\bl(\dddd:N_{k+d}\buildrel{\ddd}\over\to M_k\onto M''_k\br)=10\,(=\tau-1)\q\h{if}\q\nu_{k+5}=\dim N_{k+5}=11.
\leqno(2.4.2)$$
Here we have $\nu_{k+5}=11$ for $k\ges 8$. Since $\mu'_k=0$ for $k\ges 9$, we then get
$$\aligned\mu^{(2}_k=\nu^{(2)}_{k+5}=1&\q\h{if}\,\,\,\,k\ges 9,\\ \nu^{(2)}_{k+5}=0&\q\h{if}\,\,\,\,k\les 7,\endaligned
\leqno(2.4.3)$$
where the last assertion for $k\les 7$ follows from (2.1.1) and (1.2.4).
Note that $\mu^{(2)}_8$ and $\nu^{(2)}_{13}$ cannot be determined by these arguments as is seen by the above table.
\sk
We can show the nonvanishing (that is, the injectivity) of
$$\ddd^{(2)}:N^{(2)}_{k+10}\to M^{(2)}_k\q\h{if}\q\mu_k^{(\infty)}=0,\q\h{that is, if}\q b_f(-k/d)\ne 0.
\leqno(2.4.4)$$
Indeed, if this does not hold, then we see by using the pole order spectral sequence that there is a strictly increasing sequence $\{k_i\}_{i\in\N}$ such that $k_i-k_0\in 5\N$ and we have the nonvanishing of
$$\ddd^{(r_i)}:N^{(r_i)}_{k_i}\to M^{(r_i)}_{k_i-5r_i}\q(r_i\ges 3,\,i>0),
\leqno(2.4.5)$$
by increasing induction on $i>0$, where
$$k_i-5r_i=k_{i-1}-10,\q\h{that is},\q k_i-k_{i-1}=5(r_i-2)\q(i>0).$$
(Note that, in order to have the nonvanishing of (2.4.5), we must have the vanishing of $\ddd^{(2)}:N^{(2)}_{k_i}\to M^{(2)}_{k_i-5r_i}$.)
However, the existence of such a sequence implies that $\sum_k\mu^{(\infty)}_k$ becomes strictly smaller compared with the above table. (This can be possible since it is an infinite sequence.) This is a contradiction. So the injectivity of (2.4.4) follows.
\sk
We then get the above tables. Here the calculation in (1.4) is used for $f_1$. In fact, we cannot determine the rank of $\ddd^{(2)}:N^{(2)}_{18}\to M^{(2)}_{8}$ as in the case of $f_2$ without knowing that $-8/5$ is not a root of $b_{f_1}(s)$ by a computer calculation as in (1.4).
\msn
{\bf Remark~2.5.} In an earlier version of this paper we studied the case
$$f_1=x^5+xy^3z+y^4z+xy^4,\q\q f_2=x^5+xy^3z+y^4z.$$
Here we have
$$9/5\in R^{\,0}_{f_2}\setminus R^{\,0}_{f_1}.$$
This was first obtained by using RISA/ASIR, and was confirmed by U.~Walther using Macaulay2. 
As for the Hilbert series of the Milnor algebras, we have
$$\begin{array}{rccccccccccccccccccccccccccccc}
k:&3& 4& 5& 6& 7& 8& 9& 10& 11& 12& 13 &\cdots\\
\gamma_k:&1 &3 &6 &10 &12 &12 &10 &6 & 3& 1\\
\mu'_k:& & & & 1 &2 &2 &1 \\
\mu''_k:&1 &3 &6 &9 &10 &10 &10 &10 &10 & 10 & 10 &\cdots\\
\nu_k:& & & & & & &1 &4 &7 & 9 & 10 &\cdots\\
\mu_k:&1 &3 &6 &10 &12 &12 &11 &10 &10 & 10 & 10 &\cdots\\
\end{array}$$
The calculation of the $\mu_k$ was first done by A.~Dimca using Singular, and was confirmed by using Macaulay2 later.
In fact, this can be done quite easily by using Macaulay2 as follows:
\ms
R=QQ[x,y,z]; f=x$^{\land}$5+x*y$^{\land}$3*z+y$^{\land}$4*z+x*y$^{\land}$4;
\sk
hilbertSeries(R/(diff(x,f), diff(y,f), diff(z,f)), Order $=>$ 15)

\bs\bs
\vbox{\centerline{\bf Appendix}
\bsn
In this Appendix we give a proof of a variant of an assertion in an old version of \cite{Wa} without assuming the knowledge of logarithmic differential forms very much for nonspecialists.}
\msn
{\bf A.1.~Theorem.} {\it Let $f$ be a homogenous polynomial of $n$ variables satisfying condition~$(1)$ in the introduction, where $n\ges 2$. Then, in the notation of the introduction, we have the inequalities
$$\dim M'_k\les\dim\Gr^{n-2}_FH^{n-1}(f^{-1}(1),\C)_{{\bf e}(-k/d)}\q\h{for}\q k\in[d+1,2d],
\leqno(A.1.1)$$
where ${\bf e}(\al):=\exp(2\pi i\al)$ for $\al\in\Q$, and $F$ is the Hodge filtration.}
\msn
{\it Proof.} Set $\C[x']:=\C[x'_0,\dots,x'_n]$ with $x'_0=z$, $x'_i=x_i$ ($i\in[1,n]$). Set
$$f'=f+z^d,\q D':=\{f'=0\}\subset X':=\C^{n+1}.$$
Let $\Om^{\prime\,\ssb}$ be the complex of algebraic differential forms on $X'$. Each $\Om^{\prime\,j}$ is isomorphic to a finite direct sum of $\C[x'](-j)$ as a graded module. Here $(p)$ denotes the shift of grading by an integer $p$ in general. Define
$$A:={\rm Ker}(\df'\wedge:\Om^{\prime\,n}(-d)\to\Om^{\prime\,n+1}),$$
As in \cite{Wa}, there are exact sequences of graded $\C[x']$-modules
$$0\to A\to\Om^{\prime\,n}\buildrel{\df'\wedge}\over\too\Om^{\prime\,n+1}(d)\to(\Om^{\prime\,n+1}/\df'\wedge\Om^{\prime\,n})(d)\to 0,$$
$$0\to A(d)\buildrel{\iota}\over\to\Om^{\prime\,n}(\log D')\buildrel{\df'\wedge}\over\too\Om^{\prime\,n+1}(d)\to 0,$$
where $\iota$ is defined by
$$A(d)\ni\eta\mapsto\eta/f'\in\Om^{\prime\,n}(\log D'):=\bl\{\eta/f'\in(1/f')\Om^{\prime\,n}\,\big|\,\df'\wedge(\eta/f')\in\Om^{\prime\,n+1}(d)\br\}.$$
Recall that, for a divisor $D$ on a smooth complex algebraic variety $X$ in general, $\om\in\Om_X^j(D)$ is called {\it logarithmic} if $\ddd\om\in\Om_X^{j+1}(D)$, or equivalently, if $\ddd g\sw\om\in\Om_X^{j+1}$ for a function $g$ defining $D$ locally on $X$, see \cite{SaK}.
\sk
As in \cite{Wa}, the above exact sequences imply the graded isomorphisms
$$H_{\mm'}^0(\Om^{\prime\,n+1}/\df'\wedge\Om^{\prime\,n})(2d)=H_{\mm'}^2A(d)=H_{\mm'}^2\Om^{\prime\,n}(\log D'),
\leqno(A.1.2)$$
since $H^j_{\mm'}\C[x']=0$ for $j\les 2$, where $\mm':=(x'_0,\dots,x'_n)\subset\C[x']$, see, for instance, \cite{Ei}. (Here we use the two short exact sequences deduced from the first exact sequence.)
\sk
Set
$$\Om^{\prime\,n-1}(\log D')^{\xi}:=\bl\{\eta\in\Om^{\prime\,n-1}(\log D')\,\big|\,\xi(\eta)=0\br\},$$
where $\xi:=\msum_i\,x'_i\dd_{x'_i}$ and $\xi(\eta)$ is the contraction of $\xi$ and $\eta$. We have a short exact sequence
$$0\to(1/f')\Om^{\prime\,n+1}\buildrel{\xi}\over\too\Om^{\prime\,n}(\log D')\buildrel{\xi}\over\too\Om^{\prime\,n-1}(\log D')^{\xi}\to0.
\leqno(A.1.3)$$
In fact, let $y_0,\dots,y_n$ be the coordinates of some affine chart of the blow-up of $X'$ along the origin such that the exceptional divisor is defined by $y_0=0$. Then the restrictions of $\xi$ and $f'$ to the Zariski-open subset $\{y_0\ne 0\}$ are given respectively by
$$y_0\,\dd_{y_0},\q y_0^d\,h',$$
with $h'$ a polynomial of $y_1,\dots,y_n$. Here $y_0^d\,h'$ can be replaced with $h'$ for the definition of logarithmic forms on $\{y_0\ne 0\}$.
Then it is rather easy to show (A.1.3) on $X^{\prime*}:=X'\setminus\{0\}$ at the Zariski-sheaf level.
\sk
We now use the commutative diagram
$$\aligned0\to\,\,H^0\bl(X',\Om_{X'}^{n+1}(D')\br)\,\,\buildrel{\xi}\over\to\,\,H^0\bl(X'&,\Om_{X'}^n(\log D')\br)\,\,\buildrel{\xi}\over\to\,\,H^0\bl(X',\Om_{X'}^{n-1}(\log D')^{\xi}\br)\,\to\,0\\
\downarrow\hskip38mm&\,\,\downarrow\hskip47mm\downarrow\\
0\to H^0\bl(X^{\prime*},\Om_{X^{\prime*}}^{n+1}(D^{\prime*})\br)\buildrel{\xi}\over\to H^0\bl(X^{\prime*}&,\Om_{X'}^n(\log D^{\prime*})\br)\buildrel{\xi\,\,}\over\to H^0\bl(X^{\prime*},\Om_{X^{\prime*}}^{n-1}(\log D^{\prime*})^{\xi}\br)\to 0\endaligned$$
where $D^{\prime*}:=D'\setminus\{0\}$. The bottom row is exact since
$$H^1\bl(X^{\prime*},\Om_{X^{\prime*}}^{n+1}(D^{\prime*})\br)=0.$$
(This vanishing follows from $H^j_{\mm'}\C[x']=0$ for $j\les 2$.) We can verify that the vertical morphisms are isomorphisms by using the Hartogs-type lemma (applied to $\Om_{X'}^j$). So (A.1.3) is given by the top row of the diagram.
\sk
Set
$$Z':=\{f'=0\}\subset Y':=\PP^n.$$
The above argument implies
$$\bl(\Om^{\prime\,n-1}(\log D')^{\xi}\br)^{\sim}=\Om_{Y'}^{n-1}(\log Z'),
\leqno(A.1.4)$$
where the left-hand side is the sheaf associated with the graded module $\Om^{\prime\,n-1}(\log D')^{\xi}$.
(This is well-known to specialists, see \cite[Proposition 2.11]{DeSc}.)
\sk
Combining (A.1.2) and (A.1.3), we get the graded isomorphisms
$$\aligned H_{\mm'}^0(\Om^{\prime\,n+1}/\df'\wedge\Om^{\prime\,n})(2d)&=H_{\mm'}^2\bl(\Om^{\prime\,n-1}(\log D')^{\xi}\br)\\
&=\mopl_{k\in\Z}\,H^1\bl(Y',\Om_{Y'}^{n-1}(\log Z')(k)\br),\endaligned$$
where the last isomorphism follows from (A.1.4) together with \cite[Proposition 2.1.5]{Gr}. In particular, we get
$$H_{\mm'}^0(\Om^{\prime\,n+1}/\df'\wedge\Om^{\prime\,n})_{2d}=H^1\bl(Y',\Om_{Y'}^{n-1}(\log Z')\br).
\leqno(A.1.5)$$
This is compatible with the action of $G:=\{\la\in\C\mid\la^d=1\}$ which is defined as in \cite{Va2} by
$$z\mapsto\la z\q(\la\in G).$$
Take an embedded resolution
$$\pi:(\Yt',\Zt')\to(Y',Z'),$$
which is equivariant for the action of $G$, see \cite{AW}. (Here $\Zt'$ is {\it reduced}.) We have the isomorphisms compatible with the action of $G$
$$H_{\mm'}^0(\Om^{\prime\,n+1}/\df'\wedge\Om^{\prime\,n})_{2d}=\mopl_{k=1}^{d-1}\,H_{\mm}^0(\Om^n/\df\wedge\Om^{n-1})_{2d-k}\otimes z^k\,\ddd z/z,
\leqno(A.1.6)$$
$$\aligned H^1\bl(\Yt',\Om_{\Yt'}^{n-1}(\log\Zt')\br)&=\Gr_F^{n-1}H^n(f^{\prime\,-1}(1),\C)_1\\&=\mopl_{k=1}^{d-1}\,\Gr_F^{n-2}H^{n-1}(f^{-1}(1),\C)_{{\bf e}(k/d)}\otimes z^k\,\ddd z/z,\endaligned
\leqno(A.1.7)$$
where the action of $G$ on the last terms of (A.1.6) and (A.1.7) are induced by the action of $G$ on $z^k$. For the last isomorphism of (A.1.7), the Thom-Sebastiani type theorem \cite{ts} is used.
\sk
Combining (A.1.5--7) with Proposition~(A.2) below, we get the inequalities of dimensions (A.1.1) under the assumption~(1) if $k\in[d+1,2d-1]$.
For $k=2d$, we apply a similar argument to $f$ instead of $f'$. Here there is no action of $G$, and we may assume $n\ges 3$, since $M'_{2d}=0$ in the case $n=2$ by using the symmetry
$$\dim M'_k=\dim M'_{nd-k}\q(k\in\N),
\leqno(A.1.8)$$
see \cite{DS2}.
So the proof of Theorem~(A.1) is reduced to Proposition~(A.2) below.
\msn
{\bf Proposition~A.2.} {\it In the above notation and assumptions, we have the inequalities
$$\dim H^1\bl(Y',\Om_{Y'}^{n-1}(\log Z')\br)^{(k)}\les\dim H^1\bl(\Yt',\Om_{\Yt'}^{n-1}(\log\Zt')\br)^{(k)}\q(k\in[0,d-1]),
\leqno(A.2.1)$$
where $(*)^{(k)}$ is the subspace of $(*)$ on which the action of $\la\in G$ is given by the multiplication by $\la^k$.}
\msn
{\it Proof.} By the Leray spectral sequence
$$E_2^{p,q}=H^p\bl(Y',R^q\pi_*\Om_{\Yt'}^{n-1}(\log\Zt')\br)\Longrightarrow H^{p+q}\bl(\Yt',\Om_{\Yt'}^{n-1}(\log\Zt')\br),$$
we get a canonical injection
$$H^1\bl(Y',\pi_*\Om_{\Yt'}^{n-1}(\log\Zt')\br)\into H^1\bl(\Yt',\Om_{\Yt'}^{n-1}(\log\Zt')\br).
\leqno(A.2.2)$$
\sk
By using condition~(1) together with the Hartogs-type lemma (applied to $\Om_{Y'}^{n-1}$), we can prove a short exact sequence
$$0\to\pi_*\Om_{\Yt'}^{n-1}(\log\Zt')\to\Om_{Y'}^{n-1}(\log Z')\to\F\to 0,$$
with
$$\dim{\rm supp}\,\F=0.$$
This implies a canonical surjection
$$H^1\bl(Y',\pi_*\Om_{\Yt'}^{n-1}(\log\Zt')\br)\onto H^1\bl(Y',\Om_{Y'}^{n-1}(\log Z')\br).
\leqno(A.2.3)$$
So Proposition~(A.2) follows, since (A.2.2) and (A.2.3) are compatible with the action of $G$. This completes the proof of Theorem~(A.1).
\msn
{\bf Remarks~A.3.} (i) A canonical inclusion of vector spaces was asserted in an old version of \cite{Wa} instead of the inequality of dimensions as in (A.1.1) (without assuming condition~(1)). This seems rather difficult to prove unless one can show the bijectivity of (A.2.3).
For applications as in \cite{DiSt}, the inequality of dimensions is sufficient.
\ms
(ii) U.~Walther has recently informed us that the definition of logarithmic forms for non-reduced divisors with normal crossings in \cite{Wa} is quite different from the one in \cite{De}. For instance, if $D$ is a smooth reduced divisor on a smooth variety $X$, then the logarithmic 1-forms for the non-reduced divisor $mD$ seems to be given by
$$\Om_X^1(\log mD)=\Om_X^1(\log D)\bl((m-1)D\br),$$
where $m\ges 2$. So the canonical morphism
$$H^1\bl(X,\Om_X^1(\log D)\br)\to H^1\bl(X,\Om_X^1(\log mD)\br)$$
does not seem to be always bijective.
(This may depend on $D$ and the normal bundle $N_{D/X}$.) 
\sk
Note also that global logarithmic forms for non-normal crossing divisors are not necessarily closed (for instance, consider $\om=-(y/f)\,\ddd x+(x/f)\,\ddd y$ on $\PP^2$ for a homogeneous polynomial $f\in\C[x,y]$ of degree $d\ges 3$, where $x,y$ are the coordinates of $\C^2\subset\PP^2$), see also \cite{Wo}, etc.
\ms
(iii) Walther uses $(f-z^d)z$ (instead of $f+z^d$) in \cite{Wa}, where the Thom-Sebastiani type theorem is not needed although the calculation of the Milnor algebra becomes a little bit nontrivial.
\ms
As a corollary of Theorem~(A.1), we get the following (which does not seem to be stated explicitly in a literature although it would be well-known to specialists):
\msn
{\bf Corollary~A.4.} {\it Assume $n=3$ and
$$\h{$Z\subset\PP^2$ is a rational cuspidal curve,}
\leqno(A.4.1)$$
where ``rational cuspidal" means that $Z$ is homeomorphic to $\PP^1$. Then}
$$\dim M'_k=0\q\h{unless}\q k\in(d,2d).
\leqno(A.4.2)$$
\msn
{\it Proof.} According to \cite[Corollary 4.3]{DP}, we have
$$\dim M'_k\les\dim M'_{k+1}\q(k<3d/2),\q\dim M'_k\les\dim M'_{k-1}\q(k>3d/2).
\leqno(A.4.3)$$
Moreover, (1.2.2) holds under the hypothesis (A.4.1). (In fact, we have the nonvanishing of the restriction morphism $H^2(\PP^2,\Q)\to H^2(Z,\Q)$, since the intersection number $H\cdot Z$ is equal to $d$, where $H$ is a hyperplane of $\PP^2$.) So (A.4.2) follows from Theorem~(A.1) for $k=2d$ together with the symmetry (A.1.8).
\msn
{\bf Remarks~A.5.} (i) In \cite{DiSt} it is conjectured that under the hypothesis (A.4.1) we have
$$\dim M'_k\les 1\,\,\,\,(\forall\,k\in\N).
\leqno(A.5.1)$$
This is shown there at least if $d$ is even. Note that, in the case where condition (1.2.3) is satisfied, we have as in (1.2)
$$\dim H^2(f^{-1}(1),\C)_{{\bf e}(k/d)}=1\q(k\in[1,d-1]).
\leqno(A.5.2)$$
It follows from (A.4.3) that, for the proof of (A.5.1) in the $d$ even case, it is actually enough to show (A.5.2) only for $k=d/2$. Then (A.5.1) can be reduced to the assertion that $-1$ is not an eigenvalue of the Milnor monodromy of isolated singularities of $Z$, see \cite{DiSt}. 
\sk
In the $d$ odd case, however, the above argument only implies that (A.5.1) holds in the case where either $e^{2\pi i(d{-}1)/2d}$ or $e^{2\pi i(d{+}1)/2d}$ is not an eigenvalue of the Milnor monodromy for all the isolated singularities of $Z$. Here we use the canonical inclusion
$$H^1(f^{-1}(1),\Q)\into\mopl_{z\in{\rm Sing}\,Z}\,H^1(F_{\!h_z},\Q)^{T^d},
\leqno(A.5.3)$$
where $F_{\!h_z}$ is the Milnor fiber of a local defining function $h_z$ of $(Z,z)$ for $z\in{\rm Sing}\,Z$, and $(*)^{T^d}$ denotes the $T^d$-invariant part of $(*)$ with $T$ the Milnor monodromy and $d:=\deg f$. In fact, the inclusion (A.5.3) follows from the short exact sequence in \cite[Theorem 0.1]{DS1} by using a well-known relation between the Milnor monodromy and the vertical (or local system) monodromy of the local system $(\Hc^{n-1}\psi_f\Q_X)|_{C\setminus\{0\}}$ in the $f$ homogeneous polynomial case, where $C$ is an irreducible component of ${\rm Sing}\,f$.
\ms
(ii) There are many examples of plane rational cuspidal curves with
$$\h{$\dim M'_k=1\,$ if $\,k\in(d,2d)$, and 0 otherwise.}
\leqno(A.5.4)$$
For instance, this holds if
$$f=x^ay^b+z^d\q\h{with}\q a,b>0,\,\,\,a+b=d.$$
This easily follows from the Thom-Sebastiani type lemma for Milnor algebras, see also \cite{DiSt}. Here we need the condition $(a,b)=1$ to assure that the curve is irreducible.
\sk
However, there are many quite complicated examples of plane rational cuspidal curves, see \cite{AD}, \cite{DiSt}, etc.
For the moment, it is unclear whether we have always for such curves
$$\Sp(f)=t^{(d+1)/d}+\cdots+t^{(2d-1)/d}\,\,?
\leqno(A.5.5)$$
In case condition (1.2.3) is satisfied, this is equivalent to the question:
$$\Gr_F^pH^2(f^{-1}(1),\C)=0\q\h{unless}\q p=1\,\,?
\leqno(A.5.6)$$
Note that (A.5.5) is closely related to an improvement of the assertion in Remark~(i), see Remarks~(iii) and (iv) below.
\ms
(iii) Assume (A.5.5) is true under the assumption (A.4.1). Then the inequality (A.5.1) holds with $d$ odd if the following condition (A.5.7) is satisfied for $\al=d'/d$ or $d''/d$, where $d':=(d-1)/2$, $d'':=(d+1)/2$.
$$\h{$\al$ is not a spectral number for any singular point of $Z$.}
\leqno(A.5.7)$$
\ms
In fact, this assertion follows from Theorem~(A.1) and (A.5.5) together with (A.4.3) and (A.1.8), since (A.5.3) is a morphism of mixed Hodge structure compatible with the action of $T$.
Note that the above assumption is satisfied, for instance, if any spectral number of isolated singularities of $Z$ which is contained in $\Z\,\frac{1}{d}$ is greater than $1/2$, see Remark~(v) below for a formula of the spectral numbers in terms of Puiseux pairs.
\ms
(iv) Assume (A.5.5) is true under the assumption (A.4.1). Then the inequality (A.5.1) holds if the multiplicity of the eigenvalue $e^{2\pi i(d{-}1)/2d}$ of the Milnor monodromy is at most $1$ in the strong sense; that is, if there is at most one singular point of $Z$ such that $e^{2\pi i(d{-}1)/2d}$ is an eigenvalue of the Milnor monodromy and moreover its multiplicity is one.
In fact, it is enough to show that the assumption in Remark~(iii) is satisfied. If $(d{-}1)/2d$ is a spectral number, then $(3d{+}1)/2d$ is by the symmetry of the spectral numbers, and hence $(d{+}1)/2d$ is not by the hypothesis on the multiplicity of the eigenvalue of the monodromy. Note that the monodromy is defined over $\Q$, and $e^{2\pi i(d{+}1)/2d}$ is the complex conjugate of $e^{2\pi i(d{-}1)/2d}$ so that the same condition on the multiplicity holds.
\ms
(v) In this paper we say that $(k_1,n_1),\dots,(k_g,n_g)$ are the {\it Puiseux pairs} of an irreducible plane curve germ if it has a Puiseux expansion
$$\aligned y&=\msum_{1\les i\les[k_1/n_1]}\,a_{0,i}\,x^i\\
&+\msum_{0\les i\les[k_{2}/n_2]}\,a_{1,i}\,x^{(k_1+i)/n_1}\\
&+\msum_{0\les i\les [k_{3}/n_3]}\,a_{2,i}\,x^{k_1/n_1 +(k_2+i)/n_1n_2}\\
&\q\q\q\q\q\vdots\\
&+\msum_{i\ges 0}\,a_{g,i}\,x^{k_1/n_1+k_2/n_1n_2+\cdots+(k_g+i)/n_1{\cdots}n_g}\endaligned
\leqno(A.5.8)$$
where $a_{j,i}\in\C$ with $a_{j,0}\ne 0\,\,(j\in[1,g])$, and $k_j,n_j\in\Z_{>0}$ with ${\rm GCD}(k_j,n_j)=1$ and $n_j>1$.
Note that this definition is slightly different from the one in \cite{Le}, etc. (This is in order to simplify the inversion formula of Puiseux pairs in the proof of (A.5.9) below.)
\sk
Let $w_i$ be integers defined inductively by
$$w_1=k_1,\q w_i=w_{i-1}\,n_{i-1}\,n_i+k_i\q(i>1).$$
By \cite{exp} the Steenbrink spectral numbers which are smaller than 1 (and are counted with multiplicities) can be expressed as follows:
$$\biggl\{\biggl(r+\frac{i}{n_{\nu}}+\frac{j}{w_{\nu}}\biggr)\frac{1}{n_{\nu +1}\cdots n_g}\biggr\}
\leqno(A.5.9)$$
where $i,j,r,\nu\in\N$ with $i\in[1,n_{\nu}-1]$, $j\in[1,w_{\nu}-1]$, $r\in[0,n_{\nu +1}\cdots n_g-1]$, $\nu\in[1,g]$, and $i/n_{\nu}+j/w_{\nu}<1$. 
(This is compatible with a formula for zeta functions in terms of Puiseux pairs \cite{Le} by forgetting the integer part of the spectral numbers.)

\end{document}